\documentstyle[12pt]{article}

\def\pplogo{\vbox{\kern-\headheight\kern -17pt
\halign{##&##\hfil\cr&{IASSNS-HEP-93/43}\cr
\rule{0pt}{2.5ex}&{CLNS 93/1237}\cr
\rule{0pt}{2.5ex}&September, 1993\cr}
}}

\newcommand{\pfit}[1]{{\it #1:\/}}
\newcommand{\qedsymbol}{Q.E.D.}

\newenvironment{pf*}{\noindent\pfit }{\qedsymbol\par\par\medskip\par}

\makeatletter
\def\title{\@dblarg\@xtitle}
\def\@title{}
\def\@xtitle[#1]#2{\gdef\@title{\quad \noexpand{\vglue1pt} #2}}
\def\author{\@dblarg\@xauthor}
\def\@author{}
\def\@xauthor[#1]#2{\gdef\@author{#2}}
\def\@address{}
\def\address#1{\gdef\@address{#1}}
\def\@email{}
\def\email#1{\gdef\@email{#1}}

\def\trailer{\bigskip\par\noindent {\ixpt {\sc \@address}}}
\date{}
\def\dedicatory#1{\def\@date{\normalsize\it#1}}
\def\subjclass#1{\def\@thefnmark{}\@footnotetext{1991
    {\it Mathematics Subject Classification.} #1}}
\def\keywords#1{\def\@thefnmark{}\@footnotetext{
    {\it Key words and phrases.} #1}}
\def\thanks#1{\def\@thefnmark{}\@footnotetext{#1}}

\def\ps@firstpage{\ps@plain \def\@oddhead{\hss\pplogo}%
  \let\@evenhead\@oddhead % in case an article starts on a left-hand page
}
\def\maketitle{\par
 \begingroup
 \def\thefootnote{\fnsymbol{footnote}}
 \def\@makefnmark{\hbox
 to 0pt{$^{\@thefnmark}$\hss}}
 \if@twocolumn
 \twocolumn[\@maketitle]
 \else \newpage
 \global\@topnum\z@ \@maketitle \fi\thispagestyle{firstpage}\@thanks
 \endgroup
 \setcounter{footnote}{0}
 \let\maketitle\relax
 \let\@maketitle\relax
 \gdef\@thanks{}\gdef\@author{}\gdef\@title{}\let\thanks\relax}

\def\thebibliography#1{\section*{References\@mkboth
 {REFERENCES}{REFERENCES}}\small\list
 {\arabic{enumi}.}{\settowidth\labelwidth{[#1]}\leftmargin\labelwidth
 \advance\leftmargin\labelsep
 \usecounter{enumi}}
 \def\newblock{\hskip .11em plus .33em minus .07em}
 \sloppy\clubpenalty4000\widowpenalty4000
 \sfcode`\.=1000\relax}

\makeatother

\newif\iffn\fnfalse
\makeatletter
\@ifundefined{reset@font}{\let\reset@font\empty}{} %needed for ancient LaTeX
\long\def\@footnotetext#1{\insert\footins{\reset@font\footnotesize
    \interlinepenalty\interfootnotelinepenalty
    \splittopskip\footnotesep
    \splitmaxdepth \dp\strutbox \floatingpenalty \@MM
    \hsize\columnwidth \@parboxrestore
   \edef\@currentlabel{\csname p@footnote\endcsname\@thefnmark}\@makefntext
    {\rule{\z@}{\footnotesep}\ignorespaces
      \fntrue#1\fnfalse\strut}}}
\makeatother

\newfont{\bbbfont}{msbm10 scaled\magstep1}  % msbm12 does not exist
\newfont{\smallbbbfont}{msbm8}
\newfont{\tinybbbfont}{msbm6}
\newfont{\footbbbfont}{msbm10}
\newfont{\smallfootbbbfont}{msbm7}
\newfont{\tinyfootbbbfont}{msbm5}

\newcommand{\Bbb}[1]{\iffn
    \mathchoice{\mbox{\footbbbfont #1}}{\mbox{\footbbbfont #1}}
    {\mbox{\smallfootbbbfont #1}}{\mbox{\tinyfootbbbfont #1}}\else
    \mathchoice{\mbox{\bbbfont #1}}{\mbox{\bbbfont #1}}
    {\mbox{\smallbbbfont #1}}{\mbox{\tinybbbfont #1}}\fi}

\newcommand{\rtimes}{\mathbin{\mbox{\bbbfont\char"6F}}}
\newcommand{\semidirect}{\rtimes}

\newfont{\symfont}{msam10 scaled\magstep1}  % msam12 does not exist
\newcommand{\dashrightarrow}{\mathrel{\mbox%
   {\symfont\char"39\char"39\char"4B}}}

\newcommand{\text}[1]{\mathchoice{\mbox{\rm #1}}{\mbox{\rm #1}}
    {\mbox{\scriptsize\rm #1}}{\mbox{\tiny\rm #1}}}
\newcommand{\operatorname}[1]{\mathop{\rm #1}\nolimits}
\newcommand{\overset}{\stackrel}
\newcommand{\eqref}[1]{{\rm (\ref{#1})}}

\newcommand{\catquot}{\mathchoice
{\mathrel{\mskip-4.5mu/\!/\mskip-4.5mu}}
{\mathrel{\mskip-4.5mu/\!/\mskip-4.5mu}}
{\mathrel{\mskip-3mu/\mskip-4.5mu/\mskip-3mu}}
{\mathrel{\mskip-3mu/\mskip-4.5mu/\mskip-3mu}}
}

\newcommand{\cB}{{\Upsilon}}
\newcommand{\C}{\Bbb{C}}
\newcommand{\cD}{{\cal D}}

\newcommand{\cK}{{\cal K}}
\newcommand{\cM}{{\cal M}}
\renewcommand{\O}{{\cal O}}
\renewcommand{\P}{\Bbb P}
\newcommand{\Q}{\Bbb Q}
\newcommand{\R}{\Bbb R}
\newcommand{\Z}{\Bbb Z}
\newcommand{\Aut}{\operatorname{Aut}}
\newcommand{\Auttilde}{\mathop{\widetilde{\rm Aut}}\nolimits}
\newcommand{\Hom}{\operatorname{Hom}}
\newcommand{\Ker}{\operatorname{Ker}}
\newcommand{\Coker}{\operatorname{Coker}}
\newcommand{\suchthat}{\ |\ }

\newcommand{\gen}{\operatorname{gen}}
\renewcommand{\div}{\operatorname{div}}
\newcommand{\Div}{\operatorname{Div}}
\newcommand{\WDiv}{\operatorname{WDiv}}
\newcommand{\ad}{\operatorname{ad}}
\newcommand{\CPL}{\operatorname{CPL}}
\newcommand{\cpl}{\operatorname{cpl}}
\newcommand{\DR}{{\operatorname{DR}}}
\newcommand{\Pic}{\operatorname{Pic}}

\newcommand{\polyhedron}{P}
\newcommand{\polar}{\polyhedron^\circ}
\newcommand{\fan}{\Delta}
\newcommand{\polarfan}{\fan^\circ}
\newcommand{\normal}{{\cal N}}
\newcommand{\finiteset}{{\Xi}}
\newcommand{\toric}{_{\text{toric}}}
\newcommand{\poly}{_{\text{poly}}}
\newcommand{\Vdual}{\widehat{V}^\circ}
\newcommand{\GIT}{GIT}
\newcommand{\PL}{PL}

\newtheorem{theorem}{Theorem} 
\newtheorem{lemma}{Lemma} 
\newtheorem{definition}{Definition} 
\newtheorem{conjecture}{Conjecture} 

\begin{document}
\title{The Monomial-Divisor Mirror Map}
\author{Paul S. Aspinwall, Brian R. Greene  and David R. Morrison}
\address{Aspinwall:\ School of Natural Sciences, Institute for Advanced
Study, Princeton, NJ  08540\\
Greene:\ F.R. Newman Laboratory of Nuclear Studies, Cornell University,
Ithaca, NY 14853\\
Morrison:\ School of Mathematics, Institute for Advanced Study,
Princeton, NJ 08540}
\date{}
\thanks{Research partially supported by
 DOE grant
DE-FG02-90ER40542, a National Young Investigator award,
NSF grant DMS-9103827, and
 an American Mathematical Society Centennial Fellowship.
}

\renewcommand{\LARGE}{\Large\bf}

\maketitle

\renewcommand{\Large}{\large}

\begin{abstract}
For each family of Calabi-Yau hypersurfaces in toric varieties,
Batyrev has
proposed a possible mirror partner (which is also a family of
Calabi-Yau hypersurfaces).
We explain a  natural construction
 of the isomorphism
between certain Hodge groups of these hypersurfaces, as predicted
by mirror symmetry, which we call the {\em monomial-divisor mirror map}.
We indicate how this map can be
interpreted as the differential of the expected mirror isomorphism between
the moduli spaces of the two Calabi-Yau manifolds.
We formulate a very precise conjecture about the form of that mirror
isomorphism,
which when combined with some earlier conjectures of the third author
 would completely specify it.
We then
conclude that the moduli spaces of the nonlinear
sigma models whose targets are the different birational models of a
Calabi-Yau space should be connected by analytic continuation,
and that further analytic continuation should lead to moduli spaces
of other kinds
of conformal field theories.
(This last conclusion was first drawn by Witten.)
\end{abstract}

\section{Reflexive polyhedra}

Mirror symmetry, which proposes that Calabi-Yau manifolds should come
in pairs with certain remarkable
properties, is a phenomenon that was first observed
in the physics literature \cite{dixon,LVW,CLS,GP}.%
\footnote{For general mathematical background on
mirror symmetry and mirror pairs, we refer the reader to
\cite{mirrorbook} and \cite{guide}.}
The most concrete realization of this phenomenon---actually the only
one in which there is a physical argument linking the conformal
field theories associated to the pair of Calabi-Yau manifolds---is given
by the Greene-Plesser
orbifolding construction \cite{GP} for
Fermat hypersurfaces in weighted projective spaces
and certain quotients
of them by finite groups.  Roan \cite{roan-mirror,roan-topological}
has given a natural description of this construction in terms of
toric geometry, and he showed that the mirror phenomenon in that case
can be interpreted as a kind of duality between toric hypersurfaces.
This enabled him to give rigorous mathematical proofs of certain
formulas discovered by physicists.

Batyrev \cite{batyrev1}
has recently found an elegant characterization of Calabi-Yau
hypersurfaces which are ample Cartier divisors in (mildly singular)
toric varieties.
The characterization is stated in terms of the
{\em Newton polyhedron} of the hypersurface, which is the convex hull
of the monomials appearing in its equation.
This is always an integral polyhedron, that is,
 a compact convex polyhedron $\polyhedron$ whose vertices are elements
of a lattice $M$ in a real affine space $M_{\R}:=M\otimes\R$.
Batyrev's characterization states that
the general hypersurface with Newton polyhedron $\polyhedron$ is Calabi-Yau
(that is, has trivial canonical bundle and at worst Gorenstein
canonical singularities),
provided
 that $0$ is in the interior of $\polyhedron$, and  that each
affine hyperplane $H\subset M_{\R}$
which meets $\polyhedron$ in a face of codimension one
has the form
\[H:=\{y\in M_{\R}\suchthat\langle \ell,y\rangle=-1\}\]
for some $\ell=\ell(H)$ in the dual lattice $N:=\Hom(M,\Z)$.
An integral polyhedron with this property is called {\em reflexive}.

The normals $\ell(H)$ of supporting hyperplanes $H$ for codimension-one
faces of a reflexive polyhedron $\polyhedron$ have as their convex
hull the {\em polar polyhedron} $\polar$, which is defined to be
\[\polar:=\{x\in N_{\R}\suchthat
\langle x,y\rangle\ge-1 \text{ for all }  y\in\polyhedron\}.\]
Batyrev showed that the polar polyhedron $\polar$ of a reflexive integral
polyhedron $\polyhedron$ is itself a reflexive integral polyhedron
(with respect to the dual lattice $N$).
This led him to propose that hypersurfaces $X$ and $Y$ with Newton
polyhedra $\polyhedron$ and $\polar$, respectively, should form a
 mirror pair.

The evidence for Batyrev's proposal is of several kinds.  First and
foremost is the fact that this polar polyhedron construction specializes
to Roan's interpretation of the Greene-Plesser
orbifolding construction  in
the case of quotients of Fermat hypersurfaces in weighted projective spaces.
This is encouraging, since as noted above the
Greene-Plesser construction provides the {\em only} complete example of mirror
symmetry for hypersurfaces---the only example for which there is a physical
argument for the existence of a mirror isomorphism of the corresponding
conformal field theories.  A second piece of evidence (which we discuss
more fully below) is an isomorphism between certain  Hodge groups associated
to $X$
and $Y$ (extending the work of Roan),
as would be predicted by mirror symmetry.  And finally,
Batyrev shows that his construction is compatible with the existence
of certain ``quantum symmetries'' as expected based on physical reasoning.
This quantum
symmetry behavior looks somewhat unnatural mathematically, so verifying
it is an important check.

This evidence falls short of fully establishing a mirror symmetry
relationship between $X$ and $Y$, since it does not link the
corresponding conformal field theories.  However, it does provide
 strong
grounds for suspecting the existence of  a mirror isomorphism.  And the
 naturality of Batyrev's polar polyhedron construction
is  extremely compelling (at least to
mathematicians).

If mirror symmetry does hold between $X$ and $Y$, there will be an isomorphism
between Hodge groups $H^{1,1}(\widehat{X})$ and $H^{d-1,1}(\widehat{Y})$,
where $\widehat{X}\to X$ and $\widehat{Y}\to Y$ are appropriate (partial)
resolutions of singularities, and $d$ is the common
dimension of $X$ and  $Y$.
The existence of such an isomorphism had been
 shown quite explicitly by Roan \cite{roan-mirror,roan-topological}
in the weighted Fermat hypersurface case---%
the general case is addressed by Batyrev.
In the earlier preprint versions of  \cite{batyrev1}, Batyrev
found an equality between the dimensions of certain subspaces
$H^{1,1}\toric(\widehat{X})$ and
$H^{d-1,1}\poly(\widehat{Y})$
of the Hodge groups,
mistakenly believed to have been the entire spaces.  In the final version
of \cite{batyrev1}, he shows that the full Hodge groups are
isomorphic, following suggestions made by the present authors.  The
error in the earlier version of the paper was fortuitous, however,
as it revealed that the mirror isomorphism might be expected to preserve
those  subspaces.

In this note, we explain a very natural construction of the isomorphism
between $H^{1,1}\toric(\widehat{X})$ and
$H^{d-1,1}_{\text{poly}}(\widehat{Y})$ ,
and indicate how it can be
interpreted as the differential of the expected mirror map between
the moduli spaces (when restricted to appropriate subspaces of
those moduli spaces).
The space $H^{d-1,1}_{\text{poly}}(\widehat{Y})$
is isomorphic to the space
 of first-order polynomial deformations of
$Y$,
and can be generated by {\em monomials}; the space
$H^{1,1}\toric(\widehat{X})$ consists of that part of the second
cohomology of $\widehat{X}$
coming from the ambient toric variety, and can be generated by
toric {\em divisors}.  Our map comes from a natural one-to-one
correspondence between monomials and toric divisors
whose definition is inspired
by the constructions of Roan and Batyrev;
we have named it the {\em monomial-divisor
mirror map}.

\section{Divisors} \label{divisorsection}

Our first task is to describe the partial resolutions of singularities
we will use, and the divisors on them.
Let $\fan$ be a fan determining a toric variety
 (see \cite{Oda} or \cite{Fulton} for the definitions,
and for proofs of the facts we review below).
The support $|\fan|$ of $\fan$ is a subset of a real vector space
$N_{\R}$, and the convex cones $\sigma$ in the fan $\fan$ are rational
polyhedral cones with respect to a lattice $N$ in $N_{\R}$;
the algebraic torus which acts on the toric variety is
$T:=N\otimes\C^*$.
We let $\fan(1)$ denote the set of one-dimensional cones in $\fan$.
There is a natural {\em generator} map $\gen:\fan(1)\to N$ which assigns to
each
one-dimensional cone $\rho$ the unique generator $\gen(\rho)$
of the semigroup
 $\rho\cap N$.
Each such $\rho$ also has an associated $T$-invariant Weil divisor
$D_\rho$ in the toric variety,
which is the closure of the $T$-orbit corresponding to the cone $\rho$.

One can always describe a {\em projective} toric variety by beginning with
 a compact convex polyhedron $\polyhedron$ in
a real affine space $M_{\R}$,
integral with respect to a lattice $M$ in $M_{\R}$.  The
projective toric variety is then determined by the
{\em normal
fan} of the polyhedron $\polyhedron$; this
is the fan $\normal(\polyhedron)$ consisting of all cones
$\normal(\polyhedron,p)$ to $\polyhedron$ at $p\in\polyhedron$, where
\[\normal(\polyhedron,p):=\{x\in N_{\R}\suchthat
\langle x,p\rangle\le\langle x,y\rangle \text{ for all }  y\in\polyhedron\},\]
and where $N:=\Hom(M,\Z)$ is the dual lattice, and $N_{\R}:=N\otimes\R$.
Each proper face of the polar polyhedron $\polar$ of $\polyhedron$
is contained in a unique cone in $\normal(\polyhedron)$, which is
the cone over that face.

The toric variety $V$ determined by the fan $\normal(\polyhedron)$
is the natural one in which the hypersurfaces $X$ with Newton polyhedron
$\polyhedron$ are ample divisors.  In the case of a reflexive
integral polyhedron $\polyhedron$, the general such $X$
is an anti-canonical divisor in $V$ and will be a Calabi-Yau
variety with canonical singularities, as proved by Batyrev \cite{batyrev1}.
  We need to partially resolve
the singularities of $V$ while retaining the triviality of the canonical
bundle of the hypersurface,
getting as close as possible to a complete resolution.

To do this,  construct a blowup $\widehat{V}\to V$, determined by a fan
$\fan$ which is a subdivision of the fan $\normal(\polyhedron)$.  There will
be an induced blowup $\widehat{X}\to X$ of hypersurfaces, where $\widehat{X}$
is the proper transform of $X$ on $\widehat{V}$.  In order to maintain
the triviality of the canonical bundle of $\widehat{X}$,
 restrict the set $\fan(1)$
of one-dimensional cones as follows:  the image
$\finiteset:=\gen(\fan(1))$ of $\fan(1)$ in $N$
should lie in the set
$\polar\cap N$, where $\polar$
is the polar polyhedron of $\polyhedron$.
(We will sometimes restrict $\finiteset$ to lie in the subset
$(\polar\cap N)_0\subset{\polar\cap N}$
 consisting of those
lattice points in $\polar$
which do not lie in the interior of a codimension-one face of $\polar$.)
Oda and Park \cite{OP} (cf.\ also \cite{Stanley}) have
 shown the existence of a simplicial subdivision $\fan$
of the fan $\normal(\polyhedron)$ such that
$\gen(\fan(1))=(\polar\cap N){-}\{0\}$
(or any subset thereof),
and such that the corresponding $\widehat{V}$ is projective.
In general, there will be many such fans $\fan$.

Since the fan $\fan$ is simplicial, the toric variety $\widehat{V}$
has the structure of an {\em orbifold} (formerly called
{\em $V$-manifold} \cite{satake3}):
it can be covered by open sets of the form
$U/G_U$ where $G_U$ is a finite group acting on a manifold $U$ such that
the fixed locus of any $1\ne g\in G_U$ has real codimension at least $2$.
The open sets $U/G_U$ are used to define the notion of
{\em orbifold-smooth
differential forms}, pieced together from $G_U$-invariant smooth forms
on the open sets $U$.  Many of the theorems about the differential
geometry of smooth algebraic varieties have natural orbifold versions.
In particular, there are orbifold
de~Rham cohomology groups $H^k_\DR(\widehat{V},\R)$ isomorphic to
the real \v Cech cohomology \cite{satake3}, and
orbifold Hodge groups $H^{p,q}(\widehat{V})$ which satisfy a version
of the Dolbeault
theorem \cite{baily1}.
The general hypersurface $\widehat{X}\subset\widehat{V}$ is also an
orbifold, and has orbifold de Rham and Hodge groups of its own.

We can describe the group $\WDiv_T(\widehat{V})$ of
toric Weil divisors  on $\widehat{V}$
and their images in the Chow group $A_{n-1}(\widehat{V})$
(where $n=d+1$ is the dimension of $\widehat{V}$),
as follows (cf.\ Cox \cite{cox}).
There is a natural isomorphism\footnote{We use the notation
$\Z\langle S\rangle$ for the free abelian group on the set $S$,
and $\Z^S$ for the $\Z$-module of maps from $S$ to $\Z$,
which is naturally isomorphic to the dual lattice
$\Hom(\Z\langle S\rangle,\Z)$ of $\Z\langle S\rangle$.
The map determined by $\varphi\in\Z^S$ is denoted by $s\mapsto\varphi_s$.
}
 $\alpha:\Z^\finiteset\to\WDiv_T(\widehat{V})$
which sends the function $\varphi\in\Z^\finiteset$ to the divisor
\[\sum\varphi_aD_{\gen^{-1}(a)}.\]
Under this isomorphism, if we define an embedding
$\ad_\finiteset:M\to\Z^\finiteset$ by sending
$m\in M$ to the function $\ad_\finiteset(m)$ defined by
$\ad_\finiteset(m):a\mapsto\langle a,m\rangle$,
then
\[\div(\chi^m)=-\alpha(\ad_\finiteset(m)),\]
where $\chi^m:T\to\C^*$ is the character of $T$ associated to $m$, regarded
as a meromorphic function on $\widehat{V}$.   Thus, $M$ gives rise to
linear equivalences among toric divisors.  In fact,
there is  an exact sequence
\begin{equation} \label{exactseq}
0 \longrightarrow M \overset{\ad_\finiteset}{\longrightarrow}
\Z^{\finiteset} \overset{\bar{\alpha}}{\longrightarrow}
A_{n-1}(\widehat{V}) \longrightarrow 0 ,
\end{equation}
where $\bar{\alpha}$ denotes the composite of $\alpha$ with the projection
to the Chow group.
This is nothing other than the usual description of $A_{n-1}$ as
``divisors modulo linear equivalence'', since $\Z^\finiteset$ represents
toric divisors and $M$ represents the linear equivalences among them.

\medskip

To compute the group of toric divisors on the hypersurface
$\widehat{X}$, we use the natural restriction maps from
divisors on $\widehat{V}$ (which exists since each toric divisor
on $\widehat{V}$ meets $\widehat{X}$ in a subvariety of codimension $1$):
\[\begin{array}{ccccccccc}
0&\longrightarrow&M&\longrightarrow&\WDiv_T(\widehat{V})&\longrightarrow&
A_{n-1}(\widehat{V})&\longrightarrow&0\\
&&{\scriptstyle||}&&\downarrow&&\downarrow&&\\
0&\longrightarrow&M&\longrightarrow&\WDiv_T(\widehat{X})&\longrightarrow&
A_{d-1}(\widehat{X})&&
\end{array}.\]
This time, the toric divisors need not generate the entire Chow group;
we denote the image of $\WDiv_T(\widehat{X})$ in $A_{d-1}(\widehat{X})$
by $A_{d-1}(\widehat{X})\toric$.  Its complexification we call
the {\em toric part of $H^{1,1}$}, denoted by
$H^{1,1}\toric(\widehat{X}):=A_{d-1}(\widehat{X})\toric\otimes\C$.

The kernel of the
restriction map $\WDiv_T(\widehat{V})\to\WDiv_T(\widehat{X})$ is easy
to describe.  A divisor with trivial restriction must be supported on
divisors which are disjoint from the general hypersurface
$\widehat{X}\subset\widehat{V}$.  Since
the line bundle $\O_V(X)$ is generated by its global sections, the
general hypersurface $X\subset V$
will not meet the zero-dimensional strata of $V$ (in the
stratification by $T$-orbits).  So any divisor on $\widehat{V}$
which maps to such a stratum will be disjoint from $\widehat{X}$,
the proper transform of $X$.
Such divisors are characterized by the property that the
 corresponding point in $\finiteset$ lies in the interior of some
codimension-one face of $\polar$.
Other toric divisors on $\widehat{V}$ cannot be disjoint from
$\widehat{X}$, since they map to larger strata of $V$ which are not
disjoint from $X$.

Thus, if we let $\finiteset_0=\finiteset\cap(\polar\cap N)_0$
be the subset of $\finiteset$ consisting
of all points which do {\em not} lie in interiors of
codimension-one faces of
$\polar$, we find that $\WDiv_T(\widehat{X})\cong\Z^{\finiteset_0}$
and that
\begin{equation} \label{divisorsA}
A_{d-1}(\widehat{X})\toric\cong\Coker(\ad_{\finiteset_0})\cong
\Z^{\finiteset_0}/M.
\end{equation}
In particular, if $\finiteset\supset(\polar\cap N)_0-\{0\}$
then
$A_{d-1}(\widehat{X})\toric\cong\Z^{(\polar\cap N)_0-\{0\}}/M$,
and hence
\begin{equation} \label{divisors}
H^{1,1}\toric(\widehat{X})\cong(\Z^{(\polar\cap N)_0-\{0\}}/M)\otimes\C.
\end{equation}

\section{Monomials}

Our task in this section is to describe moduli spaces for hypersurfaces
in $\widehat{V}$.
We retain the notation of the previous section:  $\widehat{V}$ is the
 toric variety associated to a subdivision $\fan$ of the normal fan
$\normal(\polyhedron)$ of a reflexive polyhedron $\polyhedron$.
We assume that $\fan$ is simplicial, so that $\widehat{V}$ is
$\Q$-factorial; we also assume that $\widehat{V}$ is projective.

Given a hypersurface $\widehat{X}\subset\widehat{V}$, the space of
first order
deformations of complex structure of $\widehat{X}$ is isomorphic
to $H^1(\Theta_{\widehat{X}})$.  The simplest way to deform the
complex structure on $\widehat{X}$ is to perturb the equation of
the hypersurface; this leads to a subspace
$H^1(\Theta_{\widehat{X}})\poly\subset H^1(\Theta_{\widehat{X}})$
of {\em polynomial} first-order deformations.
(It is quite possible for this to be a proper subspace \cite{pdm}.)
In the case that $\widehat{X}\subset\widehat{V}$ is a Calabi-Yau
hypersurface, we can use the isomorphism $H^1(\Theta_{\widehat{X}})
\cong H^{d-1,1}(\widehat{X})$ to also specify a ``polynomial'' subspace
$H^{d-1,1}\poly(\widehat{X})\subset H^{d-1,1}(\widehat{X})$
of the corresponding Hodge group.

In principle, the moduli spaces of the hypersurfaces $\{\widehat{X}\}$
should be fairly easy to describe.  Global sections of
$\O_{\widehat{V}}(\widehat{X})$ provide equations for the hypersurfaces,
and the entire family can be described as
$\P H^0(\O_{\widehat{V}}(\widehat{X}))$.  But we need to mod out by
automorphisms of $\widehat{V}$, and this is where technical complications
arise.

Let $D$ be a toric divisor on $\widehat{V}$, and write
$D=\sum_{a\in\finiteset} d_a D_{\gen^{-1}(a)}$.  There is a natural
isomorphism between $H^0(\O(D))$ and the space $\C^{P_D\cap M}$,
where $P_D$ is the polytope
\[P_D:=\{y\in M_{\R} \suchthat \langle a,y\rangle\ge-d_a \text{ for all }
a\in\finiteset\} \]
(cf.~\cite{Fulton}).
In fact, if we identify $H^0(\O(D))$ with the space of meromorphic
functions on $\widehat{V}$ which have (at worst) poles along $D$,
then to each $m\in P_D\cap M$ we can associate the meromorphic
function $\chi^m$: it has at worst poles along $D$ thanks to the
definition of $P_D$.
In the special case $D=\sum_{\rho\in\fan(1)}D_\rho\in|-K_{\widehat{V}}|$,
the polytope
$P_{\Sigma\, D_\rho}$ coincides with the original polyhedron
$\polyhedron\subset M_{\R}$ used to describe $V$.

The automorphism group $\Aut(\widehat{V})$ of a $\Q$-factorial toric
variety $\widehat{V}$ has been described recently by Cox \cite{cox},
generalizing some results from the smooth case due to
Demazure \cite{demazure}.  Cox's description is in terms of a central
extension of $\Aut(\widehat{V})$ by a torus $G$, which
fits in an exact sequence
\begin{equation}\label{auttilde}
1\longrightarrow G\longrightarrow \Auttilde(\widehat{V})
\longrightarrow \Aut(\widehat{V})\longrightarrow1,
\end{equation}
where $G:=\Hom(A_{n-1}(\widehat{V}),\C^*)$.  The advantage of working
with $\Auttilde(\widehat{V})$ is that it acts naturally on all cohomology
groups $H^0(\O(D))$
at once.  The clearest way to see these actions is to follow Cox again
and introduce the {\em homogeneous coordinate ring}\/
$S:=\C[x_a]_{(a\in\finiteset)}$ of $\widehat{V}$.  This ring can be
(multi) graded by defining the {\em degree}\/ of the monomial
$\prod x_a^{\varphi_a}$ to be the divisor class
$[\sum \varphi_a D_{gen^{-1}(a)}]$
in $A_{n-1}(\widehat{V})$.  For a fixed divisor $D$, the set of
elements of degree $[D]$ in the homogeneous coordinate ring can be
identified with $H^0(\O(D))$:  the meromorphic function $\chi^m$ with
$m\in P_D\cap M$ corresponds to the homogeneous monomial
$x^{\div(\chi^m)+D}:=\prod x_a^{\langle a,m\rangle + d_a}$.

The torus $T:=\Hom(M,\C^*)$ which acts on $\widehat{V}$ is naturally
a subgroup of $\Aut(\widehat{V})$; the induced extension $\widetilde{T}$
of $T$ by $G$ has the form $\widetilde{T}:=\Hom(\Z^\finiteset,\C^*)$.
In fact, if we
apply the functor $\Hom(\mbox{---},\C^*)$ to the natural exact sequence
\eqref{exactseq},
we get a sequence for $\widetilde{T}$ which fits as the first row in
the commutative diagram
\[\begin{array}{ccccccccc}
1&\longrightarrow&G&\longrightarrow&\widetilde{T}&\longrightarrow&
T&\longrightarrow&1\\
&&{\scriptstyle||}&&\cap&&\cap&&\\
1&\longrightarrow&G&\longrightarrow&\Auttilde(\widehat{V})&\longrightarrow&
\Aut(\widehat{V})&\longrightarrow&1
\end{array}.\]

The grading of the homogeneous coordinate ring $S$ can also be described
in terms of the action of $G$ on $S$.  The torus $\widetilde{T}$ acts
on $S$ in a transparent way: each monomial in $S$ can be written in the
form $x^\varphi=\prod x_a^{\varphi_a}$ for some $\varphi\in\Z^{\finiteset}$,
and the action of $t\in\Hom(\Z^{\finiteset},\C^*)$ sends $x^\varphi$ to
$t(\varphi)\cdot x^\varphi$.  When we restrict this action to the subgroup
$G=\Hom(A_{n-1}(\widehat{V},\C^*))$, then for each divisor class $[D]$,
the subspace of $S$ on which
$G$ acts via the character $\gamma\mapsto\gamma([D])$ is precisely
$H^0(\O(D))\cong\bigoplus\C\cdot x^{\div(\chi^m)+D}$.

The induced action of $t\in\widetilde{T}$ on $H^0(\O(D))$ then sends
$\chi^m$ to $t(\alpha^{-1}(\div(\chi^m)+D))\cdot\chi^m$, for every
$m\in P_D\cap M$.  This action can be described in terms of the map
$\Z\langle P_D\cap M\rangle\to\Z^\finiteset$ defined by
\begin{equation}\label{action}
m\mapsto(a\mapsto\langle a,m\rangle+d_a),
\end{equation}
which induces a homomorphism of tori
$\widetilde{T}\to(\C^*)^{P_D\cap M}$
that determines the
action of $\widetilde{T}$ on $\C^{P_D \cap M}$.
Notice that the map \eqref{action} factors
as a composite of two maps
\begin{equation} \label{twomaps}
\Z\langle P_D\cap M\rangle\to M\oplus\Z\to\Z^\finiteset
\end{equation}
with the first map given by
$m\mapsto(m,1)$ and the second map given by
$(m,k)\mapsto \ad_\finiteset(m)+k\cdot\alpha^{-1}(D)$.
The corresponding homomorphism of tori factors as
\begin{equation} \label{torusfactor}
\widetilde{T}\to T\times\C^*\to(\C^*)^{P_D\cap M}.
\end{equation}
In the special case $D=\sum D_\rho$, the induced map
$\C^*\to(\C^*)^{\polyhedron\cap M}$ is simply the diagonal embedding.

The groups $\Aut(\widehat{V})$ and $\Auttilde(\widehat{V})$ are
not in general reductive.  Thus, to construct moduli spaces for
hypersurfaces\footnote{Batyrev \cite{batyrev2} has constructed moduli
spaces for {\em affine}\/ hypersurfaces, obtaining a somewhat
different space than ours if
$(\polyhedron\cap M)_0\ne(\polyhedron\cap M)$ (it even has a different
dimension).
Batyrev and Cox \cite{BC} have recently considered a construction
similar to the one described here.}
on $\widehat{V}$, we should use Fauntleroy's extension \cite{fauntleroy1}
 of Mumford's Geometric Invariant Theory (GIT) \cite{GIT}, and attempt to
construct a quotient for the
action of $\Auttilde(\widehat{V})$ on $H^0(\O(D))$.  It would be
interesting to know if
this construction of moduli spaces
can be carried out in general---Fauntleroy has carried it out in
some special cases \cite{fauntleroy2}.
We can at least obtain a birational model of the desired moduli
space by using a fairly standard result (cf.\ \cite{Rosenlicht,CDT})
which guarantees the existence of an $\Auttilde(\widehat{V})$-stable
Zariski-open set $U\subset H^0(\O(D))$ such that the geometric quotient
$U/\Auttilde(\widehat{V})$ exists.
We
indicate the birational class of such quotients
with the notation
$H^0(\O(D))\catquot \Auttilde(\widehat{V})$.

In the case of interest in this paper, $D=-K_{\widehat{V}}$.  To
study this particular moduli space,
 we take a simpler course of action, and restrict our attention
to a subspace of $H^0(\O(-K_{\widehat{V}}))$
 on which $\widetilde{T}$ acts in such
a way that the quotient exists and
has the ``expected'' dimension for the entire
moduli space.
In a wide class of examples, $\widetilde{T}$ is in fact the connected
component of the identity in $\Auttilde(\widehat{V})$, and the only
differences between our moduli space and the ``true'' moduli space
for hypersurfaces are
a remaining quotient by a finite group, and a possible ambiguity in the
choice of Zariski-open set
used in constructing the quotient.
In particular, the map from our space to the true moduli space is a
dominant map between spaces of the same dimension.
We hope that this latter property is true
in general, but postpone that question to a future investigation.

Our construction of a simplified model for the
hypersurface moduli
space relies on another result of Demazure and Cox about
$\Auttilde(\widehat{V})$.  They show that
\[\dim\Auttilde(\widehat{V})=\dim\widetilde{T}+\#(R(N,\fan)),\]
where
\begin{eqnarray*}
R(N,\fan):=\{m\in M&\suchthat&\langle \gen(\rho),m\rangle\le1
\text{ for all }\rho\in\fan(1),\\
&&
\text{with equality for a unique } \rho=\rho_m\in\fan(1)\}
\end{eqnarray*}
is the set of {\em roots} of the toric variety $\widehat{V}$ associated
to the fan $\fan$.
Note that for each root
$m\in R(N,\fan)$, we have $-m\in\polyhedron\cap M$.  In fact,
the set $-R(N,\fan)$ can be characterized
 as the subset of $\polyhedron\cap M$
consisting of lattice points which lie in the interiors of codimension-one
faces of $\polyhedron$.  We can thus decompose
\[\polyhedron\cap M=-R(N,\fan)\cup(\polyhedron\cap M)_0,\]
and write
\[\C^{\polyhedron\cap M}=\C^{-R(N,\fan)}\oplus\C^{(\polyhedron\cap M)_0}.\]
The subgroup $\widetilde{T}\subset\Auttilde(\widehat{V})$ preserves
this direct sum decomposition, so we can let $\widetilde{T}$ act on the
second factor $\C^{(\polyhedron\cap M)_0}$ alone.

Our ``simplified hypersurface
moduli space'' will be
the
\GIT\ quotient
$\C^{(\polyhedron\cap M)_0}_{\text{ss}}/\widetilde{T}$.
(We regard the action of $\widetilde{T}$ on $\C^{(\polyhedron\cap M)_0}$
as specifying a linearization of the action on
$\P(\C^{(\polyhedron\cap M)_0})$, so there is no ambiguity in
the choice of \GIT\ quotient.)
There is then a natural rational map
\begin{equation}\label{rationalmap}
\C^{(\polyhedron\cap M)_0}\catquot\widetilde{T}\dashrightarrow
\C^{\polyhedron\cap M}\catquot\Auttilde(\widehat{V})
\end{equation}
which could be refined to a regular map
\begin{equation}\label{map}
\C^{(\polyhedron\cap M)_0}_{\text{ss}}/\widetilde{T}\to
U/\Auttilde(\widehat{V})
\end{equation}
if an appropriate set of
``semistable'' points $U\subset\C^{\polyhedron\cap M}$
were available from (generalized) \GIT.

Note that $\Ker(\xi_{[-K]})$, which is a subgroup of both $\widetilde{T}$
and $\Auttilde(\widehat{V})$, acts trivially on both spaces.
By equation \eqref{torusfactor},
$\widetilde{T}/\Ker(\xi_{[-K]})\cong T\times\C^*$.
Note also that the two
quotient spaces can be expected to have the same dimension.

\begin{definition}
We say that the family
$\{\widehat{X}\}$ has the {\em dominance property} if the
natural rational map
$\C^{(\polyhedron\cap M)_0}\catquot\widetilde{T}\dashrightarrow
\C^{\polyhedron\cap M}\catquot\Auttilde(\widehat{V})$
is a dominant map between two varieties
of the same dimension.
(Note that this property is independent of the choice of quotients.)
\end{definition}
This dominance property clearly holds if $R(N,\fan)=\emptyset$;
we expect that it should hold
in general, but have not checked this.

The ``simplified hypersurface
moduli space'' parameterizes hypersurfaces with
equations of the form
\[\sum_{m\in(P\cap M)_0}c_m\chi^m=0\]
modulo the equivalences given by the action of
$\widetilde{T}/\Ker(\xi_{[-K]})\cong T\times\C^*$.
The $\C^*$ factor is diagonally embedded in $(\C^*)^{(\polyhedron\cap M)_0}$,
and so gives an overall scaling of the equation.
We can describe a Zariski-open
subset\footnote{The apparent lack of naturality in this step of
our construction---why restrict to a subset?---will be redressed
later in the paper.}
of our moduli space by restricting to equations with $c_0\ne0$, and
using the overall scaling of the equation
to set that coefficient $c_0$ equal to $1$.
Thus, the open subset
can be described as a quotient
$\C^{(\polyhedron\cap M)_0{-}\{0\}}\catquot T$
with a point $c\in\C^{(\polyhedron\cap M)_0{-}\{0\}}$
corresponding to the hypersurface with equation
\[\chi^0+\sum_{m\in(P\cap M)_0{-}\{0\}}c_m\chi^m=0.\]

Let $\cB_0=(\polyhedron\cap M)_0-\{0\}$ to simplify notation.
Here is the crucial observation for the construction of the
monomial-divisor mirror map: the action of $T=N\otimes\C^*$ on $\C^{\cB_0}$
is induced by tensoring the homomorphism $\ad_{\cB_0}:N\to \Z^{\cB_0}$
with $\C^*$.  (The explicit identification
of $\ad_{\cB_0}$ as the homomorphism
needed to specify the $T$-action
 follows immediately from
the definition of the maps in equation \eqref{twomaps},
since $\ad_{\cB_0}$ is dual to the natural map $\Z\langle{\cB_0}\rangle\to M$
induced by the inclusion ${\cB_0}\subset M$.)  In particular,
the tangent space to the simplified moduli space
$\C^{(\polyhedron\cap M)_0{-}\{0\}}\catquot T$
has the form
$(\Z^{(\polyhedron\cap M)_0{-}\{0\}}/N)\otimes\C$.

When the family $\{\widehat{X}\}$ has the dominance property
(e.g., when
$R(N,\fan)=\emptyset$), the induced rational map
$\C^{(\polyhedron\cap M)_0{-}\{0\}}\catquot T\dashrightarrow
\C^{\polyhedron\cap M}\catquot\Auttilde(\widehat{V})$
is also dominant.
In this case,
we can  describe the tangent space to the
space of polynomial deformations of a general Calabi-Yau hypersurface
 $\widehat{X}\subset\widehat{V}$
as
\begin{equation} \label{monomials}
T_{[\widehat{X}],\,\,\C^{\polyhedron\cap M}\catquot\Auttilde(\widehat{V})}=
H^{d-1,1}\poly(\widehat{X})\cong
(\Z^{(\polyhedron\cap M)_0{-}\{0\}}/N)\otimes\C,
\end{equation}
where $[\widehat{X}]$
represents the class of $\widehat{X}$ in
$\C^{\polyhedron\cap M}\catquot\Auttilde(\widehat{V})$.

We can apply these same considerations to the family of hypersurfaces
determined by the polar polyhedron $\polar$, which Batyrev has proposed
as a mirror partner for the family $\{\widehat{X}\}$.  To do this,
we need to choose a simplicial subdivision $\polarfan$ of
$\normal(\polar)$ which determines a projective toric variety $\Vdual$
and a family of hypersurfaces $\widehat{Y}\subset\Vdual$.
Replacing
 $\widehat{X}$, $\polyhedron$, $M$,  $N$ by
$\widehat{Y}$, $\polar$, $N$,  $M$, respectively, in equation
\eqref{monomials}, we find that
\begin{equation} \label{Ymonomials}
T_{[\widehat{Y}],\,\,\C^{\polar\cap N}\catquot\Auttilde(\Vdual)}=
H^{d-1,1}\poly(\widehat{Y})\cong
(\Z^{(\polar\cap N)_0{-}\{0\}}/M)\otimes\C,
\end{equation}
whenever $\{\widehat{Y}\}$ has the dominance property.

The monomial-divisor mirror map is now evident, when one compares
equations \eqref{divisors} and \eqref{Ymonomials}.
(Note that the same map
$\ad_{(\polar\cap N)_0-\{0\}}$
is used to define the embedding $M\to\Z^{(\polar\cap N)_0-\{0\}}$
in both cases.)

\begin{theorem} Let $\polyhedron$ be a reflexive polyhedron,
integral with respect
to the lattice $M$, and let $\polar$ be its polar polyhedron,
which is integral with respect to the dual lattice $N$.  Let $\fan$
and $\polarfan$ be simplicial subdivisions of $\normal(\polyhedron)$
and $\normal(\polar)$, respectively, and let $\{\widehat{X}\}$
and $\{\widehat{Y}\}$ be the corresponding families of Calabi-Yau
hypersurfaces.
Assume that
$(\polar\cap N)_0{-}\{0\}\subset \gen(\fan(1))
  \subset (\polar\cap N){-}\{0\}$.
Assume also that $R(M,\polarfan)=\emptyset$, or more
generally, simply assume that $\{\widehat{Y}\}$ has the dominance
property.  Then there is a natural isomorphism
\begin{equation}\label{mdmm}
H^{d-1,1}\poly(\widehat{Y})\overset{\cong}{\to} H^{1,1}\toric(\widehat{X})
\end{equation}
induced by equations \eqref{divisors} and \eqref{Ymonomials},
since both spaces are naturally isomorphic to
$\Coker(\ad_{(\polar\cap N)_0-\{0\}})\otimes\C$.
We call the isomorphism \eqref{mdmm} the {\em monomial-divisor mirror map}.
\end{theorem}

\section{K\"ahler cones}

Let $\widehat{V}$ be a $\Q$-factorial toric variety,
determined by a simplicial fan $\fan$, and let $\finiteset=\gen(\fan(1))$.
We can describe the group $\Div_T(\widehat{V})$ of toric
Cartier divisors on $\widehat{V}$ as follows.
In order for $D=\sum d_aD_{\gen^{-1}(a)}$ to be Cartier,
there must be a continuous piecewise linear (PL) function
$\psi_D:|\fan|\to\R$
(called the {\em support function determined by $D$})
which is linear on each cone $\sigma\in\fan$,
which takes integer values on $|\fan|\cap N$, and which satisfies
\begin{equation}\label{PL}
\psi_D(a)=- d_a\qquad\text{ for all }a\in\finiteset.
\end{equation}
Since the fan $\fan$ is simplicial, the \PL\ function $\psi_D$
is completely determined by the values specified in equation \eqref{PL}
and the fan $\fan$:
one just extends by linearity to each (simplicial) cone in $\fan$.
The integrality
condition $\psi_D(|\fan|\cap N)\subset\Z$ remains nontrivial, however.
(If we require instead that $\psi_D(|\fan|\cap N)\subset\Q$, we get
the group of $\Q$-Cartier divisors.) Since $\fan$ is simplicial,
every Weil divisor is $\Q$-Cartier, that is, $\widehat{V}$ is
$\Q$-factorial.  Put another way, there is a natural isomorphism
between $\Pic(\widehat{V})\otimes\Q$ and $A_{n-1}(\widehat{V})\otimes\Q$.

The ample Cartier divisors (or ample $\Q$-Cartier divisors) are characterized
by {\em strict convexity} of $-\psi_D$, where $\psi_D$ is the
support function determined by $D$ (cf.\ \cite{Fulton}).\footnote{Our
signs are chosen to conform to the literature as closely as possible,
while giving the term ``convexity'' its conventional meaning.}
This means the following:  given  a \PL\ function $\eta$
which is linear on each cone $\sigma\in\fan$, let $u_\sigma\in M=\Hom(N,\Z)$ be
the
linear function such that
\[\langle x,u_\sigma\rangle =\eta(x)\quad\text{ for all } x\in\sigma.\]
The convexity condition is:
\[\eta(x)\ge\langle x,u_\sigma\rangle\quad\text{ for all } x\in|\fan|;\]
the convexity if {\em strict} if the inequality is strict for all
$x\not\in\sigma$.
(In fact, it suffices to check this at the points of $\finiteset$:
for convexity, we must have
\[\eta(a)\ge\langle a,u_\sigma\rangle\quad\text{ for all } a\in\finiteset\]
with equality whenever $a\in\sigma$;
strict convexity requires a strict inequality whenever
$a\not\in\sigma$.)
The cone of real convex \PL\ functions is denoted by $\CPL(\fan)$,
following the notation of Oda and Park \cite{OP}.
If there exists a strictly convex function in $\CPL(\fan)$, the fan
$\fan$ is called {\em regular}.

There is
a related cone (introduced by
Gel'fand, Zelevinski\v\i, and Kapranov \cite{GKZ}):
\[\CPL^\sim(\fan)=\{\varphi\in\R^\finiteset\suchthat\exists\
\eta\in\CPL(\fan)
\text{ with } \varphi_a=\eta(a) \text{ for all } a\in\finiteset\}.\]
The definitions are constructed so that if $D$ is an ample $\Q$-Cartier
divisor, then $\alpha^{-1}(D)\in\CPL^\sim(\fan)$.  (One can choose
$\eta=-\psi_D$ in the definition.)
Note that $M_{\R}$ is contained in $\CPL^\sim(\fan)$, and the corresponding
$\eta$'s are precisely the
{\em smooth} \PL\ functions.
The image of $\CPL^\sim(\fan)$ in $\R^\finiteset/M_{\R}$,
which we may think of as the set of convex \PL\ functions modulo
smooth \PL\ functions,
 is denoted
by $\cpl(\fan)$. Under the isomorphism
$\R^\finiteset/M_{\R}\cong A_{n-1}(\widehat{V})\otimes\R
(\cong\Pic(\widehat{V})\otimes\R)$,
this cone
$\cpl(\fan)$ maps to the closed real cone generated by the ample divisor
classes on $\widehat{V}$.
An effective method of calculating all possible cones $\cpl(\fan)$
in terms of the ``linear Gale transform'' is described in \cite{OP}.

The exponential sheaf sequence gives rise to an isomorphism
$\Pic(\widehat{V})\cong H^2(\widehat{V},\R)\cong H^{1,1}(\widehat{V},\R)$,
since $h^i(\O_{\widehat{V}})=0$
for $i>0$.
Now there is a natural notion of {\em positivity} for orbifold-smooth
$(1,1)$-forms:  one requires that the $G_U$-invariant $(1,1)$-forms
on the local uniformizing sets $U$ be positive.  The K\"ahler form of
every orbifold-K\"ahler metric is easily seen to be
a positive, orbifold-smooth
$(1,1)$-form.  The set
$\cK(\widehat{V})\subset H^{1,1}(\widehat{V},\R)$
consisting of orbifold de~Rham
classes which have such a positive representative is called the
{\em K\"ahler cone} of $\widehat{V}$.

We could not find the following lemma in the literature
(although it should be known).
\begin{lemma}
Under the natural map $\Pic(\widehat{V})\to H^2_\DR(\widehat{V},\R)$
which assigns to a line bundle the corresponding orbifold de Rham class,
the ample line bundles map to positive de Rham classes.
\end{lemma}
Note that the lemma is not as obvious as is the analogous lemma
 in the smooth case, since
even if we are given a projective embedding $\widehat{V}\to\P^N$,
the pullback of the Fubini-Study form (which establishes the positivity
of the
de Rham class of $\O_{\P^N}(1)$) is {\em not}
necessarily positive
as an orbifold $2$-form.

\medskip

\begin{pf*}{Sketch of proof}
We modify an argument of Guillemin and Sternberg
\cite{GS1}.
By a result of Delzant \cite{delzant} and Audin \cite{audin},
the toric variety $\widehat{V}$ can be described
as a symplectic reduction of the action of $G=\Hom(A_{n-1}(\widehat{V}),\C^*)$
on $\C^{\finiteset}$; the specific symplectic reduction which produces
$\widehat{V}$ is $\Phi^{-1}(\alpha)/G_{\R}$, where $\Phi$ is the moment
map for the action, $G_{\R}$ is the maximal compact subgroup of $G$,
and $\alpha\in\cpl(\fan)$.  An ample line bundle $L$ on $\widehat{V}$
corresponds to a character $\chi^L$ of $G$, since $A_{n-1}(\widehat{V})$
is the character lattice of $G$; moreover, the corresponding point in
$A_{n-1}(\widehat{V})\otimes\R$ lies in $\cpl(\fan)$.  If we apply
the constructions described on pp.~520--521 of \cite{GS1}, we produce
a line bundle on $\widehat{V}$ with a specific (orbifold-)connection,
whose curvature is the symplectic form obtained by symplectic reduction
from the standard form on $\C^{\finiteset}$.  That curvature form
provides the desired positive orbifold $2$-form.
 \end{pf*}

The converse statement---that if the class of a line bundle is represented
by a positive, orbifold-smooth $2$-form then the line bundle is ample---%
is a theorem of Baily \cite{baily2}.  Putting the two together, we conclude
that the image of the cone $\cpl(\fan)$ in $H^2_\DR(\widehat{V},\R)$
is precisely the closure of the K\"ahler cone $\overline{\cK(\widehat{V})}$.

\medskip

As remarked earlier, the hypersurface $\widehat{X}$ is itself an
orbifold, so it has an orbifold K\"ahler cone
$\cK(\widehat{X})\subset H^{1,1}(\widehat{X},\R)$.
The positive orbifold-smooth K\"ahler forms
on $\widehat{V}$ corresponding to classes in the interior of
$\cpl(\fan)$ will
restrict to positive orbifold-smooth forms on $\widehat{X}$, since
$\widehat{X}$ meets all singular strata of $\widehat{V}$ transversally.
We call the resulting cone
 $\cK\toric(\widehat{X})
\subset\cK(\widehat{X})\cap H^{1,1}\toric(\widehat{X})$
 the {\em cone of toric K\"ahler classes} on $\widehat{X}$.

\section{The K\"ahler moduli space}

Mirror symmetry predicts a close relationship between the moduli space
of complex structures on one Calabi-Yau manifold $\widehat{Y}$
and the so-called ``K\"ahler moduli space'' of its mirror partner
$\widehat{X}$.  This K\"ahler moduli space, which arises in the study of
nonlinear sigma models with
target $\widehat{X}$,\footnote{The physics of these models
is believed to be as well-behaved on orbifolds as on manifolds
\cite{DHVW,DFMS}.} is
an open subset of $\cD /\Gamma$, where
\[\cD :=
\{B+i\,J\in H^2(\widehat{X},\C)\suchthat J\in\cK(\widehat{X})\},\]
and where
\[\Gamma:= H^2(\widehat{X},\Z)\semidirect \Aut(\widehat{X})\]
(cf.~\cite{compact}.)  The precise open subset of $\cD /\Gamma$
which constitutes the moduli space is difficult to determine,
since general convergence criteria for the sigma model are unknown at
present.\footnote{But as we will observe below, the region of convergence
in this particular case can be inferred from mirror symmetry.}
However, the open set is expected
to include points sufficiently
far out along any path in $\cD $ whose imaginary part is moving
towards infinity while staying away from the boundary of the K\"ahler cone.
Such paths should approach a common
point, called the {\em large radius limit},
in an appropriate partial compactification of $\cD /\Gamma$.
A general discussion of conditions under which such a limit point exists
can be found in \cite{compact}.

We will be interested in
a ``toric subspace'' of the K\"ahler moduli space,
defined by intersecting the moduli space with
$\cD\toric(\widehat{X})/\Gamma\toric$, where
\[\cD\toric(\widehat{X}):= \{B+i\,J\in H^{1,1}\toric
(\widehat{X})\suchthat J\in\cK\toric(\widehat{X})\},\]
and  $\Gamma\toric:=A_{d-1}(\widehat{X})\toric
\semidirect \Aut(\widehat{X})\toric$.
(The toric automorphisms $\Aut(\widehat{X})\toric$
are those automorphisms of $\widehat{X}$ induced by an automorphism
of the ambient toric variety $\widehat{V}$.)
If the \GIT\
of the family of hypersurfaces is well-behaved (which will be the
case for the families of primary interest to us), then for the general
hypersurface $\widehat{X}$ the group $\Aut(\widehat{X})\toric$
will be finite.

Let $L=A_{d-1}(\widehat{X})\toric$ and consider the torus
$L\otimes\C^*$ which contains $\cD/L$ as an open subset.
The cone $\cK\toric\subset L\otimes\R$, which is a rational polyhedral
cone, determines an affine torus embedding $\cM$ with a unique
$0$-dimensional orbit $p\in\cM$ under the action of the torus $L\otimes\C^*$.
We let $\overline{\cD/L}$ be the closure of $\cD/L$ in $\cM$, and
let $(\cD/L)^-$ be the interior of $\overline{\cD/L}$.
(This contains the point $p$.)  If $\Aut(\widehat{X})\toric$ is finite, then
since everything in sight is $\Aut(\widehat{X})\toric$-equivariant,
we may take the quotient and get a partial compactification $(\cD/\Gamma)^-$
of $\cD/\Gamma$ with a distinguished boundary point (again denoted by
$p$), the {\em large radius limit} point.\footnote{This construction
is identical to the one described in \cite{compact}, since Looijenga's
semi-toric compactification \cite{Looijenga} coincides with the
``toroidal embeddings'' of \cite{AMRT} when $\cK$ is rational polyhedral
and $\Gamma/L$ is finite.}
This point is the common limit of the paths described earlier.

The torus $L\otimes\C^*=A_{d-1}(\widehat{X})\toric\otimes\C^*$ which
is being compactified can be described in the form
\begin{eqnarray*}
A_{d-1}(\widehat{X})\toric\otimes\C^*&=&
(\Z^{\finiteset_0}/M)\otimes\C^*\\&=&
(\C^*)^{\finiteset_0}/(M\otimes\C^*)
\end{eqnarray*}
using equation \eqref{divisorsA}.
Now the orbits of $M\otimes\C^*$ on $(\C^*)^{\finiteset_0}$ are all good
orbits of the same dimension.
The action of $M\otimes\C^*$  on the larger space $\C^{\finiteset_0}$
may be less well-behaved, but in any case we can regard
 $(\C^*)^{\finiteset_0}/(M\otimes\C^*)$ as a representative of the
birational class of quotients
 $\C^{\finiteset_0}\catquot(M\otimes\C^*)$.

Suppose that the  family $\{\widehat{Y}\subset\Vdual\}$
associated to the polar polyhedron
of $\{\widehat{X}\}$
has the dominance property (introduced in section 3), and that $\fan$ is chosen
so that
$(\polar\cap N)_0{-}\{0\}\subset \gen(\fan(1))
  \subset (\polar\cap N){-}\{0\}$.
Then we deduce from the monomial-divisor mirror map a diagram
in which the vertical maps are dominant:
\[\begin{array}{ccccc}
\cD\toric(\widehat{X})/A_{d-1}(\widehat{X})\toric&\subset&
(\C^*)^{(\polar\cap N)_0{-}\{0\}}/(M\otimes\C^*)&=&
\C^{(\polar\cap N)_0{-}\{0\}}\catquot(M\otimes\C^*)\\
\raise6pt\hbox{$\Big\downarrow$}&&
\raise6pt\hbox{$\Big\downarrow$}&&
\overset{\overset{\scriptscriptstyle|}{\scriptscriptstyle|}}
{\scriptscriptstyle\downarrow}\\
\cD\toric(\widehat{X})/\Gamma\toric&\subset&
(A_{d-1}(\widehat{X})\otimes\C^*)/\Aut(\widehat{X})\toric&&
\C^{\polar\cap N}\catquot\Auttilde(\Vdual).
\end{array}\]
The left and center vertical maps are simply the quotient maps
by $\Aut(\widehat{X})\toric$, and the right map is the dominant rational
map from the simplified moduli space to the actual moduli space
of the family $\{\widehat{Y}\}$.

We can now formulate a ``mirror symmetry'' conjecture for these families,
which generalizes some (less precise) earlier conjectures of
Aspinwall and L\"utken \cite{AL} and Batyrev \cite{batyrev2}.

\begin{conjecture}
Suppose that $\Aut(\widehat{X})\toric$ is finite, and that
the dominance property holds for $\{\widehat{Y}\}$.  Then
\begin{enumerate}
\item
there is an open set $U\subset(\cD\toric(\widehat{X})/\Gamma\toric)^-$
containing the large radius limit point $p$ such that
$U\cap(\cD\toric(\widehat{X})/\Gamma\toric)$ is the toric part of the K\"ahler
moduli space,
\item
there is an appropriate  quotient
$\C^{\polar\cap N}\catquot\Auttilde(\Vdual)$ and a
{\em ``mirror map''}\footnote{We denote this map by $\mu^{-1}$
in order to match the conventions established in \cite{compact}.}
\[\mu^{-1}:U\to
\C^{\polar\cap N}\catquot\Auttilde(\Vdual),\]
which is an isomorphism onto its image, and which,
when restricted to $U\cap(\cD\toric(\widehat{X})/\Gamma\toric)$,
serves to identify points whose conformal field theories are
mirror-isomorphic, such that
\item
the differential of the inverse map $\mu$ at the ``large complex structure
 limit point'' $\mu^{-1}(p)$
\[d\mu:
T_{\mu^{-1}(p),\,\,\C^{\polar\cap N}\catquot\Auttilde(\Vdual)}
\to
T_{p,\,\,U}
\]
coincides with the monomial-divisor mirror map
\[
H^{d-1,1}\poly(\widehat{Y})
\to
H^{1,1}\toric(\widehat{X})
\]
up to signs,
once we have made the canonical identifications
\begin{eqnarray*}
T_{\mu^{-1}(p),\,\,\C^{\polar\cap N}\catquot\Auttilde(\Vdual)}&=&
H^{d-1,1}\poly(\widehat{Y})\\
T_{p,\,\,U}&=&H^{1,1}\toric(\widehat{X}).
\end{eqnarray*}
That is, there is some element
\[   \theta_{\fan} \in A_{d-1}(\widehat{X})\toric \otimes \C^*
                     \subset \Aut ( H^{1,1}\toric (\widehat{X}) )     \]
of order $2$, which when composed with the monomial-divisor mirror map
yields $d\mu$.  (When $d \ge 3$, the automorphism $\theta_{\fan}$ which
specifies the signs is unique.)
\end{enumerate}
In particular, the
 location of the large complex structure limit point $\mu^{-1}(p)$ can be
calculated using the monomial-divisor mirror map and the knowledge
of the cone $\cK\toric(\widehat{X})$.
\end{conjecture}

In \cite{compact}, a very general
conjecture was formulated which specifies the
``canonical coordinates'' to be used in the mirror map, up to some
constants of integration.  Those constants can be determined if one knows
the differential of the mirror map at the large radius limit point---%
for these toric hypersurfaces, that differential is supplied by the
conjecture above.
So the two conjectures together completely specify the canonical
coordinates.
A similar conjecture about canonical coordinates
for toric hypersurfaces has been independently made by Batyrev and
van Straten \cite{BvS}.

If the conjecture stated above is true, then among other things the so-called
``$3$-point functions'' (part of the conformal field theories)
must coincide under this mapping.  The  $3$-point
function on the moduli space $\C^{\polar\cap N}\catquot\Auttilde(\Vdual)$
can be calculated in terms of the variation of Hodge structure
of the family $\{\widehat{Y}\}$ (cf.~\cite{guide}),
and this gives a method to identify precisely which subset of
$\cD\toric(\widehat{X})/\Gamma\toric$ constitutes the toric part of the
K\"ahler
moduli space.  For that subset can be characterized as the domain of
convergence of a power series expansion of the $3$-point functions,
when calculated in the canonical coordinates on
$\cD\toric(\widehat{X})/\Gamma\toric$.
Using the mirror map, this power series calculation can actually be made
in $\C^{\polar\cap N}\catquot\Auttilde(\Vdual)$, by calculating
 the variation of
Hodge structure of the family $\{\widehat{Y}\}$.

The leading term in the power series expansion of the $3$-point function
on $\cD\toric(\widehat{X})/\Gamma\toric$ is the cubic form
\[H^{1,1}\toric(\widehat{X})\times H^{1,1}\toric(\widehat{X})\times
H^{1,1}\toric(\widehat{X})\to\C\]
given by the cup product; higher terms are given by
``quantum corrections'' that depend
on the point in $\cD\toric(\widehat{X})/\Gamma\toric$, and that all vanish
at the large radius limit point $p\in(\cD\toric(\widehat{X})/\Gamma\toric)^-$.
One consequence of our conjecture would therefore be an agreement between
the leading term in the variation of Hodge structure calculation for
the family $\{\widehat{Y}\}$ near the large complex structure limit point,
and the cup product cubic form on $\widehat{X}$.
This consequence was
first checked in an example by Aspinwall, L\"utken, and Ross \cite{ALR}
several years ago.\footnote{It is possible to verify that the large
complex structure limit point used in \cite{ALR} agrees with the
one predicted by the monomial-divisor mirror map.}
  More recently, Batyrev \cite{batyrev2} checked his version of this
statement in the case that $X$ itself is smooth,
and the authors \cite{AGM} checked it
in an example in which
there are five different birational choices for $\widehat{X}$ (with the same
$X\subset V$).  After learning of the results of \cite{AGM} and of the
present paper, Batyrev \cite{batyrev3} checked this consequence in general.

The $3$-point function coming from variation of Hodge structure can be
used to determine the choice of signs $\theta_{\fan}$:
the poles in that function should occur at {\em positive}\/ real values
of the canonical coordinates (cf.~\cite{AGMiii}).
The automorphism $\theta_{\fan}$ with this property can be
calculated explicitly using methods of Gel'fand, Zelevinski\v\i, and
Kapranov \cite{GKZ}; we will discuss this in detail elsewhere.

Another consequence of our conjecture is that the K\"ahler moduli
spaces for different birational models $\widehat{X}$ of the function
field of $X$ can naturally be regarded as analytic continuations of
one another.  (For after applying mirror symmetry, they are seen to
occupy different regions in the same moduli space.)  This was the
principal conclusion of our earlier paper \cite{AGM};
a similar idea is
due independently to Manin \cite{manin}.

\section{Phases and the secondary fan}

In the course of defining the
 monomial-divisor mirror map,
we made a somewhat unnatural restriction to a Zariski-open subset
of the simplified hypersurface moduli space.  We now return to
the study of the full moduli space.
The ``simplified hypersurface
moduli space'' will be birational to the quotient
$(\C^*){}^{(\polyhedron\cap M)_0}/\widetilde{T}$.
In fact, the moduli space of primary interest is the space parameterizing
those hypersurfaces whose singularities are no worse than generic.
This is the complement of the ``principal discriminant'' of
Gel'fand, Zelevinski\v\i, and Kapranov \cite{GKZ}.
One natural compactification
of the moduli space would be the one in which this ``principal discriminant''
is an ample divisor.

Whatever compactification we use,
the compactified moduli space is itself a toric variety (since it contains
the torus $(\C^*){}^{(\polyhedron\cap M)_0}/\widetilde{T}$
as a dense open subset).
If we compactify so that the principal discriminant is ample,
then the toric variety is determined by
 the Newton
polyhedron for the principal discriminant which,
as Gel'fand, Zelevinski\v\i, and Kapranov show,  has a convenient
combinatorial description as a so-called ``secondary polytope''.

To explain the combinatorics, we note that the
action of $\widetilde{T}$
on $(\C^*){}^{(\polyhedron\cap M)_0}$
 is induced by a homomorphism
$\ad^+_{(\polyhedron\cap M)_0}:N^+\to\Z^{(\polyhedron\cap M)_0}$,
dual to the map $\Z\langle{(\polyhedron\cap M)_0}\rangle\to M^+$
given  by equation \eqref{twomaps}, where $M^+=M\oplus\Z$.
We should imagine embedding the set ${(\polyhedron\cap M)_0}$
into $M^+$ via the map
$b\mapsto (b,1)$; the image is a finite set of points in the affine
hyperplane $\{(m,1)\}\subset M^+$.  The convex cone spanned by
these points we denote by $\polyhedron^+$; it is the cone over the
image of the original polyhedron $\polyhedron$ (generated by the
points ${(\polyhedron\cap M)_0}$)
in the affine hyperplane $\{(m,1)\}\subset M^+$.
The dual cone to $\polyhedron^+$ has the form $(\polar)^+$,
where $\polar$ is the polar polyhedron of $\polyhedron$.

The {\em secondary fan} (which
is the normal fan of the secondary polytope) is the fan consisting of
all cones $\cpl(\Sigma)$, where $\Sigma$ is a regular
refinement of the fan $\normal(\,(\polar)^+)$
(cf.~\cite{BFS}).
Among the possible regular refinements $\Sigma$ we find fans of the form
\[(\polarfan)^+:=
\{\text{cone over }(\sigma\cap\polyhedron)\suchthat\sigma\in\polarfan\},\]
for regular fans $\polarfan$ which refine $\normal(\polar)$. (But there
are others, which do not have this form.)  For such fans, it
is easy to see that $\cpl(\,(\polarfan)^+)=\cpl(\polarfan)$, regarding both
as cones in the same space
$\Z^{{(\polyhedron\cap M)_0}}/N^+\cong\Z^{{(\polyhedron\cap M)_0}-\{0\}}/N$.

So our chosen
compactification of the moduli space is the toric variety which is
specified by the secondary fan.  It has the pleasant property that
it includes all of the partial compactifications that were needed to
describe the ``large complex structure limits'' coming from mirror
symmetry of sigma models (since those were given by the cones
$\cpl(\polarfan)$).  It has another nice property as well, first
observed by Kapranov et al.\ \cite{KSZ} and Batyrev \cite{batyrev2}:
the compactification constructed in this way dominates all possible
  \GIT\ compactifications, coming from different choices of linearization.

What does this structure correspond to under mirror symmetry?  The
embedding $(\polar\cap N)_0\subset N$, together with a (regular)
refinement $\fan$ of
the fan $\normal(\polyhedron)$, was used to determine the projective
toric variety
$\widehat{V}$.
The
new extended embedding $(\polar\cap N)_0\subset N^+$ (whose image lies in
the affine hyperplane $\{(n,1)\}\subset N^+$), together with a regular
refinement
$\Sigma$ of the fan $\normal(\polyhedron^+)$, can also be used to
determine a toric variety, of dimension one larger than the previous variety.
Among these toric varieties we find the total spaces of canonical
bundles over the various choices of $\widehat{V}$ (when we take $\Sigma$
of the form $\fan^+$).

This is precisely the structure that Witten has found to be relevant
in his study of Landau-Ginzburg theories and their deformations
\cite{phases}.  Each choice of fan $\Sigma$ determines a different
``phase'' of the physical theory.  When the fan $\Sigma$ is of the form
$\fan^+$, the physical theory is related to the nonlinear sigma model
with target $\widehat{X}$ (where $\widehat{X}\subset\widehat{V}$
is generic).  But for other fans $\Sigma$, the physical theory is quite
different.  (We refer the reader to \cite{phases} for more details.)
Thus, not only can the different sigma models with birationally equivalent
targets be viewed as analytic continuations of one another, there are
further analytic continuations (to regions in the moduli space corresponding
to $\cpl(\Sigma)$ for $\Sigma\ne\fan^+$) to other kinds of physical
theory.

\ifx\undefined\bysame
\newcommand{\bysame}{\leavevmode\hbox to3em{\hrulefill}\,}
\fi

\trailer

\end{document}